 
\catcode`\@=11
\def\undefine#1{\let#1\undefined}
\def\newsymbol#1#2#3#4#5{\let\next@\relax
 \ifnum#2=\@ne\let\next@\msafam@\else
 \ifnum#2=\tw@\let\next@\msbfam@\fi\fi
 \mathchardef#1="#3\next@#4#5}
\def\mathhexbox@#1#2#3{\relax
 \ifmmode\mathpalette{}{\m@th\mathchar"#1#2#3}%
 \else\leavevmode\hbox{$\m@th\mathchar"#1#2#3$}\fi}
\def\hexnumber@#1{\ifcase#1 0\or 1\or 2\or 3\or 4\or 5\or 6\or 7\or 8\or
 9\or A\or B\or C\or D\or E\or F\fi}

\newdimen\ex@
\ex@.2326ex
\def\varinjlim{\mathop{\vtop{\ialign{##\crcr
 \hfil\rm lim\hfil\crcr\noalign{\nointerlineskip}\rightarrowfill\crcr
 \noalign{\nointerlineskip\kern-\ex@}\crcr}}}}
\def\varprojlim{\mathop{\vtop{\ialign{##\crcr
 \hfil\rm lim\hfil\crcr\noalign{\nointerlineskip}\leftarrowfill\crcr
 \noalign{\nointerlineskip\kern-\ex@}\crcr}}}}
\def\varliminf{\mathop{\underline{\vrule height\z@ depth.2exwidth\z@
 \hbox{\rm lim}}}}

\font\tenmsa=msam10
\font\sevenmsa=msam7
\font\fivemsa=msam5
\newfam\msafam
\textfont\msafam=\tenmsa
\scriptfont\msafam=\sevenmsa
\scriptscriptfont\msafam=\fivemsa
\edef\msafam@{\hexnumber@\msafam}
\mathchardef\dabar@"0\msafam@39
\def\dashrightarrow{\mathrel{\dabar@\dabar@\mathchar"0\msafam@4B}}
\def\dashleftarrow{\mathrel{\mathchar"0\msafam@4C\dabar@\dabar@}}

\font\tenmsb=msbm10
\font\sevenmsb=msbm7
\font\fivemsb=msbm5
\newfam\msbfam
\textfont\msbfam=\tenmsb
\scriptfont\msbfam=\sevenmsb
\scriptscriptfont\msbfam=\fivemsb
\edef\msbfam@{\hexnumber@\msbfam}
\def\Bbb#1{{\fam\msbfam\relax#1}}
\def\widehat#1{\setbox\z@\hbox{$\m@th#1$}%
 \ifdim\wd\z@>\tw@ em\mathaccent"0\msbfam@5B{#1}%
 \else\mathaccent"0362{#1}\fi}
\font\teneufm=eufm10
\font\seveneufm=eufm7
\font\fiveeufm=eufm5
\newfam\eufmfam
\textfont\eufmfam=\teneufm
\scriptfont\eufmfam=\seveneufm
\scriptscriptfont\eufmfam=\fiveeufm

\newsymbol\subsetneq 2328
\newsymbol\complement 107B
\newsymbol\ndiv 232D

\catcode`\@=12

\magnification=\magstep1
\font\title = cmr10 scaled \magstep2

\font\smalltext = cmr7
\font\smallmath= cmmi7
\font\tinymath=cmmi5
\font\smallsym = cmsy7
\font\author = cmcsc10

\font\byabs = cmr7

\parindent=1em
\baselineskip 15pt

\newcount\refcount
\newcount\seccount
\newcount\sscount
\newcount\eqcount
\newcount\boxcount
\newcount\testcount
\newcount\bibcount
\boxcount = 128
\seccount = -1

\def\sec#1{\advance\seccount by 1\bigskip\goodbreak\noindent
	{\bf\number\seccount.\ #1}\medskip \sscount = 0\eqcount = 0}
\def\subsec{\advance\sscount by 1\medskip\goodbreak\noindent
	{{\number\seccount.\number\sscount}\ \ } \eqcount = 0}
\def\namess#1{\advance\sscount by 1\global
	\edef#1{\number\seccount.\number\sscount}
	\medskip\goodbreak\noindent{\bf #1\ \ } \eqcount = 0}
\def\proc#1#2{\advance\sscount by 1\eqcount = 0
	\medskip\goodbreak\noindent{\author #1}
	{\tenrm{\number\seccount.\number\sscount}}:\ \ {\it #2}}
\def\nproc#1#2#3{\advance\sscount by 1\eqcount = 0\global
	\edef#1{#2\ \number\seccount.\number\sscount}	
	\medskip\goodbreak\noindent{\author #2}
	{\tenrm{\number\seccount.\number\sscount}}:\ \ {\it #3}}
\def\proof{\medskip\noindent{\it Proof:\ \ }}

\def\eql#1{\global\advance\eqcount by 1\global
	\edef#1{(\number\seccount.\number\sscount.\number\eqcount)}\leqno{#1}}
\def\ref#1#2{\advance\refcount by 1\global
	\edef#1{[\number\refcount]}\setbox\boxcount=
	\vbox{\item{[\number\refcount]}#2}\advance\boxcount by 1}
\def\biblio{{\frenchspacing
	\bigskip\goodbreak\centerline{\bf REFERENCES}\medskip
	\bibcount = 128\loop\ifnum\testcount < \refcount
	\goodbreak\advance\testcount by 1\box\bibcount
	\advance\bibcount by 1\vskip 4pt\repeat\medskip}}
\def\tightmatrix#1{\null\,\vcenter{\normalbaselines\mathsurround=0pt
	\ialign{\hfil$##$\hfil&&\ \hfil$##$\hfil\crcr
	\mathstrut\crcr\noalign{\kern-\baselineskip}
	#1\crcr\mathstrut\crcr\noalign{\kern-\baselineskip}}}\,}

\def\emph{\bf}
\def\colon{{:}\;}
\def\|{|\;}


\def\C{{\Bbb C}}
\def\R{{\Bbb R}}
\def\Z{{\Bbb Z}}
\def\Q{{\Bbb Q}}
\def\F{{\Bbb F}}
\def\R{{\Bbb R}}
\def\N{{\Bbb N}}

\def\LP{\left(}

\def\cF{{\cal F}}
\def\End{{\rm End}}

\def\GL{{\rm GL}}

\def\Stab{{\rm Stab}}
\def\Norm{{\rm Norm}}

\def\qed{\hfill\nobreak\rlap{$\sqcap$}$\sqcup$}

\def\ddownarrow{\rlap{$\downarrow$}\raise 2pt\hbox{$\downarrow$}}

\ref\De{P. Deligne: La Conjecture de Weil, II, {\it Inst. Hautes \'Etudes 
	Sci. Publ. Math.} {\bf 52} (1980), 138--252.}
\ref\DS{P.~Diaconis, M.~Shahshahani: On the eigenvalues of random matrices. 
	Studies in applied probability.  {\it J. Appl. Probab.}  {\bf 31A}  (1994), 49--62.}
\ref\Katz{N.~Katz: {\it Gauss sums, Kloosterman sums, and monodromy groups.}
	Annals of Mathematics Studies, 116. Princeton University Press, Princeton, NJ, 1988.}
\ref\LP{M.~Larsen, R.~Pink: Determining representations from invariant 
	dimensions, {\it Invent. Math.} {\bf 102} (1990) 377--398.}
\ref\Margulis{G.~A.~Margulis: Discrete subgroups of semisimple Lie groups. Ergebnisse der 	Mathematik und ihrer Grenzgebiete (3), 17. Springer-Verlag, Berlin, 1991.}
\ref\Richardson{Richardson, R.~W., Jr.: Deformations of Lie subgroups and the variation
	of isotropy subgroups.  {\it Acta. Math.} {\bf 129} (1972), 35--73.}
\ref\Weyl{H.~Weyl: {\it The Classical Groups}, Princeton University Press, 
	Princeton, 1939.}  

\centerline{\title Rigidity in the Invariant Theory}
\centerline{\title of Compact Groups}
\bigskip
\centerline{ Michael Larsen\footnote*
{Partially supported by NSF Grant No. DMS-0100537}}
\centerline{ Department of Mathematics}
\centerline{Indiana University}
\centerline{Bloomington, IN 47405 USA}
\bigskip
\centerline{\byabs ABSTRACT}
\smallskip
{\byabs \narrower\narrower
\textfont0 = \smalltext
\textfont1 = \smallmath
\scriptfont1 = \tinymath
\textfont2 = \smallsym
A compact Lie group $G$ and a faithful complex representation $V$ determine
the Sato-Tate measure $\mu_{G,V}$ on {\sevenmsb C}, defined as the direct image of
Haar measure with respect to the character map $g\mapsto {\rm tr}(g|V)$.
We give a necessary and sufficient condition for a Sato-Tate measure to be an isolated point in the set of all Sato-Tate measures, regarded as a subset of the space of distributions on {\sevenmsb C}.  In particular we prove that if $G$ is connected and semisimple and $V$ is irreducible, 
then $\mu_{G,V}$ is an isolated point.\par
}
\bigskip
\sec{Introduction}
Given a compact Lie group $G$ and a finite dimensional, faithful, complex 
representation $V$ of $G$, there is an associated {\emph Sato-Tate
measure} $\mu_{G,V}$ on the complex plane defined as the direct image 
of Haar measure with respect to the 
map $G\to\C$ given by the trace of $V$.  This measure is compactly supported 
and is therefore uniquely 
determined by its moments
$$\int_{\C}z^a \bar z^b\mu_{G,V}=\dim(V^{\otimes a}\otimes{V^*}^{\otimes b})^G,
\quad a,b\ge0.$$
To what extent is the pair $(G,V)$ determined by 
the measure $\mu_{G,V}$ or equivalently by the dimensions of spaces of 
$G$-invariants of $V^{\otimes a}\otimes{V^*}^{\otimes b}$? Certainly
$(G,V)$ is not always determined up to isomorphism.
For example, when $G$ is a finite group and $V$ is its regular 
representation, the only information encoded in $\mu_{G,V}$ is the order 
of $G$.  In \LP, we showed that $(G,V)$ is not in general determined up to 
isomorphism by $\mu_{G,V}$ even when $G$ is connected and semisimple.    
If, however, in addition $V$ is irreducible, then the pair is uniquely determined by its
measure.

This paper does not directly address the general problem of 
when $\mu_{G,V}$ determines $(G,V)$.  
It considers a related question: for what pairs 
$(G,V)$ is $\mu_{G,V}$ isolated from all other Sato-Tate measures in the 
space of distributions on $\C$?  Equivalently, when is the 
{\emph Sato-Tate function} 

$$F_{G,V}\colon
(a,b)\mapsto \dim(V^{\otimes a}\otimes{V^*}^{\otimes b})^G$$  
isolated in the space of all Sato-Tate functions in $\N^{\N\times\N}$ (endowed with the 
Tychonoff topology)?  It may be worth pointing out that not every limit of a 
convergent sequence of Sato-Tate measures is again a Sato-Tate measure.
For example, if $(G_i,V_i)$ is a sequence of classical groups in the same series
and  $V_i$ is the natural representation of $G_i$, then $\mu_{G_i,V_i}$ converges
to a complex (resp. real) Gaussian distribution if the $G_i$ are unitary
(resp. orthogonal or symplectic) \DS.

The original motivation for this problem comes from the cohomological theory of exponential
sums (see, e.g., \Katz~Chap.~3).   Suppose that $U$ is a geometrically connected
curve over a finite field $\F$, $\bar U$ the curve obtained
from $U$ by extension of scalars to $\overline{\F}$, and $\pi_1^a$ and $\pi_1^g$
are the fundamental groups of $U$ and $\bar U$ respectively.    
Let $\cF$ be a  lisse $\ell$-adic sheaf
on $U$, i.e., a continuous representation $\pi_1^a\to\GL_n(\Q_\ell)$.
For simplicity, assume that $\cF$ is pure of weight zero, and that $\F$ is
large enough that the Zariski closures $G^a$ and $G^g$ of $\pi_1^a$ and $\pi_1^g$ respectively
coincide.  Then, by \De,   $G^a = G^g$ has semisimple identity component.
Let $G$ denote a compact real form of the complexification of $G^a$
(with respect to a fixed embedding $\Q_\ell\hookrightarrow\C$), so $\cF$ determines
a complex faithful $n$-dimensional representation $V$ of $G$.   Every closed
point $u$ in $U$ determines a conjugacy class in $G$ whose trace in $V$ coincides
with the trace of the Frobenius in $\cF$.   Deligne's geometric analogue of the Sato-Tate
conjecture \De\ asserts that the distribution of these conjugacy classes is uniform with respect to
Haar measure on $G$.  At the trace level, this means that the distribution of traces of
Frobenius with respect to $\cF$ is given by the Sato-Tate measure $\mu_{G,V}$.
In fact, Deligne proves a strong version of this estimate, which allows one to compute
any given moment $F_{G,V}(a,b)$ effectively.   Ideally, we would like to determine $(G,V)$, but failing that, we would like to know that a finite amount of computation is enough to determine $\mu_{G,V}$.  The main theorem of the paper is insufficient for this more limited goal in general, but it turns out to be enough when $G$ is connected and $V$ is irreducible.

The situation for monodromy of $\ell$-adic sheaves over function fields is special because 
the Riemann hypothesis is a theorem, and therefore we can achieve certainty concerning individual moments of $\mu_{G,V}$.  Suppose we have an unknown compact group $G$, arising in nature as a subgroup of $\GL_n$ (where $n$ may be unknown) for which we can sample elements of $G$ uniformly distributed with respect to Haar measure and observe the trace.  Repeated experiments will lead to greater and greater certainty concerning any given moment of $\mu_{G,V}$.  We would like to know when they will lead to greater and greater certainty concerning the actual distribution $\mu_{G,V}$.

To understand how a Sato-Tate measure could fail to be isolated,  we note that from the standpoint of invariant theory, it is hard to distinguish $U(1)$ from $\Z/n\Z$ when $n$ is large.
On the other hand,  Jordan's 
structure theorem for finite subgroups of $\GL_n(\C)$ asserts, in 
effect, that a compact Lie group can be well approximated by finite 
subgroups {\it only} if its identity component is a torus.  More generally we 
expect that a compact Lie group will be well approximated by proper closed 
subgroups exactly if its identity component fails to be semisimple.  This phenomenon
turns out to explain all cases in  which a measure $\mu_{G,V}$ is a limit of other Sato-Tate measures.

The paper is organized as follows.  The first section gives two basic 
finiteness results for conjugacy classes of closed subgroups of compact Lie groups.
If not known, they are at least variants of well-known statements.  The second section proves that
$\mu_{G,V}$ cannot be isolated if
the center of $G^\circ$ has positive dimension. 
The proof of the main theorem is in the third section.  
The fourth section examines the question of when $\mu_{G,V}$ is a 
limit of distributions associated to representations of finite groups.  

I would like to thank N.~Katz for introducing me to this circle of ideas in 
his 1985--86 course at Princeton.  This paper is motivated by some more recent
questions of his. 

\sec{Some finiteness theorems}
Throughout this paper, $G$ will denote a compact Lie group, $G^\circ$
its identity component, $D$ the derived group of $G^\circ$, $Z$ the 
center of $G^\circ$, and $Z^\circ$ its identity component.  
Thus, $D$ is semisimple, 
$Z^\circ$ is a torus, and $G^\circ=DZ^\circ$. 

We begin with a version of Richardson's rigidity theorem \Richardson.

\nproc\RR{Lemma}{Let $X$ be a variety over $\R$, $G$ a linear algebraic group over
$\R$, and $H$ an $\R$-subgroup of $G$.  Let $J\subset G\times X$ denote
a closed $\R$-subvariety which is an $X$-subgroup scheme (of the constant group 
scheme $G\times X / X$.  Suppose that $H(\R)$ is compact.  Then there exists
a finite collection $S$ of closed subgroups of $H(\R)$ such that for every $x\in X(\R)$,
the intersection of $H(\R)$ with $J_x(\R)$ is conjugate in $H(\R)$ to a group in $S$.}

\proof
Let $K = J\cap H\times X$.  By a standard Noetherian induction argument, 
there is a partition of $X$ into finitely many locally closed strata over 
each of which $K$ is a smooth group scheme.  Without loss of generality,
therefore, we may assume $K$ is smooth over $X$.  Therefore, $K_x(\R)$ is a real
analytic family of subgroups of $H(\R)$ in the sense of \Richardson.
Now, $K_x(\R)$ is compact for all $x$, so every representation of $K_x(\R)$ is
completely reducible, and this means $H^1(K_x(\R), V) = 0$ for any representation
$V$.  By \Richardson~Theorem~3.1, for every $x_0\in X(\R)$, there exists 
a neighborhood $U$ such that $K_x(\R)$ is conjugate to $K_{x_0}(\R)$ 
in $H(\R)$ for all $x\in U$.  As the real locus $X(\R)$ has finitely many components in
the strong topology, there exists a finite set of compact subgroups of $H(\R)$
representing all $K_x(\R)$ up to conjugation.

\qed

Now consider a compact Lie group $G$, a faithful complex representation
$V$ of $G$, and a finite collection $\Sigma$ of pairs 
$(m,n)\in\N^2$.  If $H\subset G$ is a closed subgroup,  $V$ can 
also be regarded as a representation of $H$, and
$$F_{G,V}(m,n)\le F_{H,V}(m,n)\eqno{(\ast)}$$
for all $(m,n)$, in particular those belonging to $\Sigma$.  

\nproc\InvariantFiniteness
{Proposition}{Given $(G,V)$ and $\Sigma$ as above, there exists a 
 set $S$ of subgroups $K\subset G$ consisting of finitely many 
 $G$-conjugacy classes and such that for any 
subgroup $H\subset G$, either $F_{G,V}$ coincides with $F_{H,V}$ 
for all $(m,n)\in\Sigma$ or $H$ is contained in some $K\in S$.}

\proof For each $(m,n)\in \Sigma$, we consider the Grassmannian $X_{m,n}$ of 
$F_{G,V}(m,n)+1$-planes in $V_{m,n}:=V^{\otimes m}\otimes {V^*}^{\otimes n}$.  If, 
indeed, the inequality ($\ast$) is strict, then the space of $H$-invariants 
includes some plane indexed by $X_{m,n}$.  Let $X$ denote the disjoint 
union over $\Sigma$ of $X_{m,n}$, and for $x\in X_{m,n}$, let $K_x$ denote 
the intersection in $\GL(V_{m,n})$ of $G$ and $\Stab(x)$.  By \RR,
the resulting possibilities for $K_x$ are finite, up to conjugation in $G$.
\qed 

\medskip 
The last result in this section is a generalization to compact Lie groups of Jordan's 
theorem on finite subgroups of $GL_n(\C)$.  

\nproc\Jordan{Theorem}{If $G$ is a compact Lie group for which $G^\circ$ 
is semisimple, then there exists a finite set $S$ of proper compact 
subgroups of $G$ 
such that every proper compact subgroup can be conjugated into a subgroup 
of an element of $S$.}

\proof Any subgroup which does not meet every component of $G$ is contained
in a proper subgroup of $G$ which is a union of components of $G$.  
Therefore, without loss of generality, we may consider only subgroups $H$ 
of $G$
which meet every component of $G$. 
As $G^\circ$ is semisimple, it has finitely many normal subgroups, and if 
$H\cap G^\circ$ contains a positive-dimensional normal subgroup of $G^\circ$, 
it contains a positive-dimensional normal subgroup of $G$.    
If a maximal proper subgroup $H$ of $G$ contains a normal subgroup $K$ 
of $G$, then $H/K$ is a maximal proper subgroup of $G/K$.  By induction on 
dimension, we may consider only subgroups $H$ of $G$ which do not contain 
a positive dimensional normal subgroup of $G$.  

If $K$ is any positive dimensional characteristic subgroup of $H$, 

$$H\subset \Norm_G K\subsetneq G,$$
so without loss of generality, we may assume $H=\Norm_G K$.  If $H^\circ$ 
is noncommutative, we take $K$ to be the derived group of $H^\circ$; if 
$H^\circ$ is a non-trivial torus, we set $K=H^\circ$.  

In the first case $H$ is the normalizer of a semisimple subgroup of 
$G^\circ$, so it suffices to prove there are only finitely many such 
subgroups up to conjugacy.  There are finitely many possible isomorphism 
classes for $K$.  Given a fixed embedding $G\hookrightarrow U(n)$, 
each inclusion $K\hookrightarrow G$ determines an $n$-dimensional representation of 
$K$.  There are finitely many isomorphism classes of $n$-dimensional 
representations of a fixed semisimple group; each one determines a 
family of subgroups $K_x\subset U(n)$ indexed by 
$U(n)$ itself.  By \RR, 
there are only finitely many possibilities for $K_x \cap G\subset G$ and therefore
finitely many possibilities for $K=K_x$ and for $H$, up to 
$G$-conjugacy.

If $K$ is a torus, $H$ contains the centralizer of the torus and therefore 
some maximal torus $T$ of $G^\circ$.   Thus $K$ is a maximal torus of
$G^\circ$, which means this case contributes only one conjugacy class
of maximal subgroups.

All that remains is the case that $H$ is finite.  Up to 
conjugation, there are finitely many conjugacy classes of finite subgroups 
of $G$ of any given order: 
there are finitely many isomorphism classes of finite groups of a given 
order, and for each fixed isomorphism class, by \RR, there are finitely many 
conjugacy classes of homomorphisms to $G$.  
By Jordan's theorem, there exists a 
constant $J$ depending only on $G$ such that 
$H$ contains a normal abelian subgroup $A$ of 
index $\le J$.  If $|H| > J|Z|[G:G^\circ]$, then $A$ must contain a
non-central element $a\in G^\circ$.  The centralizer of $a$ in $G$ is
a proper closed subgroup of positive dimension, and it contains $H$, which is
absurd.   Thus, $|H|$ is bounded above, and the theorem follows.\qed

\sec{Approximating tori by torsion subgroups}

Let $Z[n]$ denote the kernel of multiplication by $n$ in $Z$.  This is a 
characteristic subgroup of $G$.  We regard the $Z[n]$ as finite subgroups
approximating $Z$.  
We would like to show that whenever $Z$ is infinite, $G$ itself can be well 
approximated by proper closed subgroups.

\nproc\TorusTorsion{Lemma}{For every positive integer $n$ there exists a 
compact subgroup $G_n$ of $G$ such that 
$G_n Z=G$, 
and $G_n\cap Z$ contains $Z[n]$ as a subgroup of finite index.  In 
particular, if $Z$ is infinite, $G_n$ is a proper subgroup of $G$.}

\def\cont{{}}

\proof Suppose that $H_\cont^2(G/Z^\circ,Z^\circ)$ is finite.  Let $m$ annihilate this
group.  Then the cohomology sequence associated to the short exact sequence
$$0\to Z^\circ[m]\to Z^\circ\to Z^\circ\to 0$$
implies every class of $H_\cont^2(G/Z^\circ,Z^\circ)$, in particular the one given by the
extension $0\to Z^\circ\to G\to G/Z^\circ\to 0$, lies in the image of
$H_\cont^2(G/Z^\circ,Z^\circ[m])$.  It  can therefore be trivialized by pulling back to
a finite central extension
$F$ of $G/Z^\circ$.  This means 
that the pull-back of the homomorphism $G\to G/Z^\circ$ to $F$ splits, 
and the splitting gives a homomorphism
$\phi\colon F\to G$ such that 
$\phi(F)Z^\circ=G$.  We can then take $G_n=\phi(F) Z[n]$.  

To prove the desired finiteness, first consider the short exact sequence
$$0\to X_*(Z^\circ)\to X_*(Z^\circ)\otimes\R\to Z^\circ\to 0,$$    
where $X_*$ denotes the group of cocharacters.
By the long exact cohomology sequence, it suffices to prove that
$H_\cont^q(G/Z^\circ,X_*(Z^\circ))$ is finite for $q=2,3$  or equivalently,
that $H_\cont^q(G/Z^\circ,\Z)$ is finite for the same values of $q$.
The first is obvious. The second is a consequence
of Whitehead's lemma.
\qed

\nproc\Nonisolation{Proposition}
{If $Z$ is infinite, then $\mu_{G,V}$ is a limit point of the sequence 
$$\{\mu_{G_n,V}\}_{n=1,2,3,\ldots}.$$
}

\proof Clearly 
$$F_{G,V}(a,b)\le F_{G_n,V}(a,b)$$
for all $a,b\in\N$.  It suffices to show that for any fixed pair $(a,b)$, 
equality holds for all but finitely many $n\in\N$.  Let $Z^\circ=
Z_1\times\cdots\times Z_k$, where each $Z_i$ is isomorphic to $U(1)$.
Let $X_{a,b}$ denote the
Grassmannian of $F_{G,V}(a,b)+1$-planes in 
$V_{a,b}:=V^{\otimes a}\otimes{V^*}^{\otimes b}$.  For $x\in X_{a,b}$, 
we let $Z_{i,x}$ denote the intersection in $\GL(V_{a,b})$ of 
$Z_i$ and $\Stab(x)$.
As these intersections are the fibers of a morphism of finite type, there is 
a uniform upper bound $N$ such that if $Z_{i,x}$ has more than $N$ points,
it must be all of $Z_i$.  For $n> N$, $G_n\subset\Stab(x)$ implies
$Z_i[n]\subset\Stab(x)$ and therefore $Z_i\subset\Stab(x)$ for all $i$.
As $G=G_n Z^\circ$, this implies $G\subset\Stab(x)$ contrary to hypothesis.  
Thus, $F_{G,V}(a,b)=F_{G_n,V}(a,b)$.\qed

\sec{Isolated measures}

In this section, we prove the main theorem of this paper, namely the 
isolation of any $F_{G,V}$ which is not a limit point by 
virtue of \Nonisolation. 

\nproc\crude{Lemma}{For any finite-dimensional Hermitian vector space $V$,
there exists $N$ such that for $a>N$,
$$F_{U(V),V}(a,a)>(\dim V - 1)^{2a}.$$
}

\proof
\def\tr{{\rm tr}}
Let $du$ denote Haar measure on $U(V)$ and
$$K=\Bigl\{u\in U(V)\Bigm\vert |\tr(u)|>\dim V-{1\over 2}\Bigr\}.$$
By the continuity of the trace function, $K$ has positive measure with 
respect to $du$, so if $a\gg 0$,
$$\eqalign{F_{U(V),V}(a,a)&=\int_{U(V)} |\tr(u)|^{2a} du
\ge\int_K |\tr(u)|^{2a} du\cr
&\ge \left(\int_K du\right)
\left(\dim V-{1\over 2}\right)^{2a}>(\dim V-1)^{2a}.}
$$
\qed

\proc{Proposition}{For all pairs $(G,V)$ there exists a neighborhood $U$ 
of $F_{G,V}$ in the Tychonoff topology such that for all $(G',V')$ with 
$F_{G',V'}\in U$, $\dim V=\dim V'$.}

\proof We have for all $a$ the inequalities 

$$F_{U(V),V}(a,a)\le F_{G,V}(a,a)\le (\dim V)^{2a}.$$
By the second main theorem of invariant theory \Weyl~2.14.A, if $a\le \dim V$, 

$$F_{U(V),V}(a,a)=a!>(a/e)^a.$$
Therefore, if $\dim V'>e(\dim V)^2$, setting $a=\lceil e(\dim V)^2\rceil$, 

$$F_{G,V}(a,a)\le(\dim V)^{2a}\le (a/e)^a<F_{U(V'),V'}(a,a)\le F_{G',V'}(a,a).$$ 
This gives a neighborhood of $F_{G,V}$ which contains only $F_{G',V'}$ 
with $\dim V'$ bounded above.  For any fixed value $\dim V'\neq \dim V$, 
by \crude, there exists $a$ such that 
$$F_{G,V}(a,a)\le(\dim V)^{2a}<F_{U(V'),V'}(a,a)\le F_{G',V'}(a,a).$$
or 
$$F_{G',V'}(a,a)\le(\dim V')^{2a}<F_{U(V),V}(a,a)\le F_{G,V}(a,a)$$  

Thus, in some neighboorhood $U$ of $F_{G,V}$, all Sato-Tate sequences 
arise from representations of degree $\dim V$.\qed 

\medskip
We can now prove the main theorem of this paper.  

\nproc\mainThm{Theorem}
{Let $\mu$ denote any Sato-Tate measure for which 
$\mu=\mu_{G,V}$ implies that $G^\circ$ is semisimple.  Then $\mu$ is an 
isolated point in the space of Sato-Tate measures.  Equivalently, there 
exists $N$ such that if 
$$F_{G,V}(a,b)=\int_{\C}z^a\bar z^b\mu$$
whenever $a,b\le N$, this equation holds for all $a$ and $b$.}

\proof Let $\mu=\mu_{G_0,V_0}$ and $n=\dim V_0$.  By the preceding proposition, all 
$\mu_{G,V}$ sufficiently close to $\mu$ satisfy $\dim V=n$. Consider the 
set of all isomorphism classes of pairs $(G,V)$ with $\dim V=n$.  We 
construct a directed graph on this set as follows: for each conjugacy 
class of inclusion maps $G_1\hookrightarrow G_2$ such that the restriction 
of the representation $V_2$ to $G_1$ is isomorphic to $V_1$, there is an 
arrow from $(G_1,V_1)$ to $(G_2,V_2)$.  If $G_2^\circ$ is semisimple and
$G_1$ is a maximal 
proper closed subgroup, we color the arrow red.  If $G_2^\circ$ is not 
semisimple, by hypothesis, there exists a  minimal non-negative integer $N$
such that 
$$F_{G_2,V_2}(a,b)\neq\int_{\C}z^a\bar z^b\mu\quad$$
for some $a+b = N$.
Let $\Sigma=\{(a,b)\in\N^2\mid a+b\le N\}$.  Consider the set $S$ of subgroups
of $G_2$ given by \InvariantFiniteness.  If $(G_1,V_1)$ is isomorphic to 
$(K,V_2)$ for some $K\in S$, then we color the arrow blue.  All other 
arrows remain uncolored.   

We consider the maximal connected subgraph $\Gamma$ 
containing the pair $(U(V),V)$ for an $n$-dimensional vector space and 
consisting only of red and blue arrows.  The vertices of this graph are
those in the original graph which can be obtained from $(U(V),V)$ by
following a sequence of colored arrows.  We 
claim this graph is finite.  It suffices to show that there are no infinite 
sequences of arrows arranged head to tail and that each pair is the head 
of finitely many arrows.  The first point is clear since any decreasing chain 
of nested  compact groups stabilizes.  The second point follows from 
\InvariantFiniteness\ and \Jordan.  Note that  $\Gamma$ contains no directed cycle since 
every  injective endomorphism of a finite-dimensional Lie group is an isomorphism.  

Assume $\mu$ is a limit point of Sato-Tate measures.  It is then a limit 
point of Sato-Tate measures associated to $n$-dimensional representations.  
For any $n$-dimensional pair $(G,V)$ in $\Gamma$
it makes sense to ask whether $\mu$ 
is a limit point of Sato-Tate measures {\emph subordinate to $(G,V)$}, that is,
of the form $\mu_{G',V}$, where 
$G'$ ranges over subgroups of $G$.  For the pair $(U(V),V)$, the answer is 
affirmitive.  Let $(G,V)$ denote a minimal pair in $\Gamma$ for which the answer 
remains affirmitive.  

If $G^\circ$ is semisimple, then any infinite sequence of Sato-Tate measures 
subordinate to $(G,V)$ contains an infinite subsequence subordinate to 
$(G_1,V)$ for some maximal subgroup $G_1$ of $G$. This contradicts
the minimality of $(G,V)$. Thus, $G^\circ$ is not semisimple.  Choose $N$ so that
for some $a,b\le N$, 
$$F_{G,V}(a,b)\neq \int_{\C}z^a\bar z^b\mu.$$
If $F_{G,V}(a,b)$ exceeds $\int_{\C}z^a\bar z^b\mu$, the same will be true for
any restriction of $V$ to a subgroup of $G$, contrary to the assumption that
$\mu$ is a limit of Sato-Tate measures subordinate to $(G,V)$.  Thus, 
$$F_{G,V}(a,b)<\int_{\C}z^a\bar z^b\mu,$$
and every sequence of Sato-Tate measures subordinate to $(G,V)$ and converging 
to $\mu$ contains a
subsequence $\mu_{G_1,V_1}$, $V_1=V\vert_{G_1}$,
such that $F_{G_1,V_1}(a,b)>F_{G,V}(a,b)$.
Every term in the subsequence is subordinate to a pair which is the target
of a blue arrow from $(G,V)$; by the minimality of $(G,V)$, this gives a 
contradiction.  \qed 

\proc{Corollary}{If $G$ is a semisimple group and $V$ is an irreducible 
representation, then $\mu_{G,V}$ is isolated. }

\proof It suffices to prove that if $\mu_{G',V'}=\mu_{G,V}$, then $G'$ has 
a semisimple identity component.  

As 
$$\dim\End_{G'}(V')=\dim(V'\otimes{V'}^*)^{G'}=F_{G',V'}(1,1)=F_{G,V}(1,1)
=\dim\End_G(V)=1,$$
the representation $V'$ is irreducible.  Consider the restriction of $V'$ 
to the identity component of $G'$.  This decomposes as a direct sum of 
isotypic factors $V'_i$, and the component group $G'/{G'}^\circ$ acts 
transitively on the set of factors.  It is a theorem of Jordan that for 
every non-trivial transitive group action, there is a group element 
without fixed points.  This element corresponds to a connected component 
of $G'$ on which the trace of $V'$ is zero.  Thus, the measure 
$\mu_{G',V'}$ has an atom at $0$.  This is impossible, since the pre-image 
of $0$ under the character of $V$ is a proper subvariety of the connected 
real-algebraic variety $G$, and therefore has measure zero \Margulis~I~2.5.3~(i).  
We conclude  that the restriction of $V'$ is isotypic of type $W'$.  As $V'$ is a 
faithful representation of $G'$, the center of ${G'}^\circ$ lies in the 
center of $G'$.  If it has positive dimension, then $F_{G',V'}(a,b)$ is 
zero, except when $a=b$.  This is impossible, since $F_{G,V}(a,b)$ is 
positive whenever $a$ and $b$ are sufficiently large and $a-b$ annihilates 
the center of $G$.  The corollary follows.\qed

\sec{Limits of finite groups}

\nproc\finiteSub{Lemma}
{If $G$ is a compact Lie group such that $D$ is non-trivial, there exists
a finite set $S'$ of proper closed subgroups such that every finite subgroup 
of $G$ is conjugate to a subgroup of an element of $S'$.}

\proof 
Let $S$ be the finite set of subgroups of $G/Z^\circ$ obtained
by applying \Jordan\ to that group.  Let $S'$ denote the set of preimages
with respect to $G\to G/Z^\circ$ of elements of $S$.  The image of any finite
subgroup is a proper closed subgroup of $G/Z^\circ$ and therefore a subgroup
of some conjugate of an element of $S$; it follows that any finite subgroup
of $G$ can be conjugated into a subgroup of an element of $S'$.
\qed

\proc{Theorem}
{If $\mu$ is a Sato-Tate measure which is the limit a a sequence of
Sato-Tate measures of representations of finite groups, then $\mu=\mu_{G,V}$
for some group $G$ such that $G^\circ$ is a torus.}

\proof
We follow the proof of \mainThm, using \finiteSub\ instead of \Jordan.
That is, for every $(G,V)$ with $D$ non-trivial, we draw a red
arrow to each $(G',V)$, $G'\in S'$.  For every $(G,V)$
with $G^\circ$ a torus, we find $a+b=N$ such that
$F_{G,V}(a,b)\neq\int_{\C}z^a\bar z^b\mu$ and draw a blue arrow to
each $(G',V)$ obtained by applying \InvariantFiniteness\ to
$\Sigma=\{(a,b)\in\N^2\mid a+b\le N\}$.  We choose a minimal $(G,V)$
on the resulting finite graph for which there exists an infinite sequence
of pairs $(G_i,V_i)$ with $G_i$ finite, subordinate to $(G,V)$ and with
Sato-Tate measures converging to $\mu$.  The contradiction results as
before.
\qed

We remark that, although it seems likely that when $\mu_{G,V}=\mu_{G',V'}$
either both or neither of $G^\circ$ and ${G'}^\circ$ are tori, this has
not been proven, so in principle for some non-toric $(G,V)$ we could still
have $\mu_{G,V}$ the limit of Sato-Tate measures of finite groups.

\biblio

\end